%% file: franklrevised.tex
\documentclass{amsart}
\usepackage{amssymb,amsmath,color}
\usepackage{dsfont}
\usepackage{mathrsfs}
\usepackage{enumitem}

\usepackage{color}
\usepackage{tikz}
\usepackage{pgf}
\usepackage{bm}
\usepackage{empheq}
\usepackage{lscape}

\newtheorem{theorem}{Theorem}[section]
\newtheorem{conjecture}[theorem]{Conjecture}
\newtheorem{lemma}[theorem]{Lemma}

\newtheorem{proposition}[theorem]{Proposition}
\theoremstyle{definition}
\newtheorem{definition}[theorem]{Definition}

\theoremstyle{remark}

\newcommand{\boxedeq}[1]{\begin{empheq}[box={\fboxsep=6pt\fbox}]{align*}#1\end{empheq}}

\date{\today}
\title{New Conjectures For Union-Closed Families}

\author{Jonad Pulaj}
\address{Konrad-Zuse-Zentrum,
Takusstra\ss e 7,
Berlin 14195, Germany}
\email{pulaj@zib.de}

\author{Annie Raymond}
\address{Department of Mathematics, University of Washington, Box
  354350, Seattle, WA 98195, USA} \email{raymonda@uw.edu}

\author{Dirk Theis}
\address{Institute of Computer Science, University of Tartu, Juhan Liivi 2, 50409 Tartu, Estonia} \email{dirk.theis@ut.ee}

\begin{document}

\begin{abstract} The Frankl conjecture, also known as the union-closed sets conjecture, states that in any finite non-empty union-closed family, there exists an element in at least half of the sets. From an optimization point of view, one could instead prove that $2a$ is an upper bound to the number of sets in a union-closed family on a ground set of $n$ elements where each element is in at most $a$ sets for all $a,n\in \mathbb{N}^+$. Similarly, one could prove that the minimum number of sets containing the most frequent element in a (non-empty) union-closed family with $m$ sets and $n$ elements is at least $\frac{m}{2}$ for any $m,n\in \mathbb{N}^+$. Formulating these problems as integer programs, we observe that the optimal values we computed do not vary with $n$. We formalize these observations as conjectures, and show that they are not equivalent to the Frankl conjecture while still having wide-reaching implications if proven true. Finally, we prove special cases of the new conjectures and discuss possible approaches to solve them completely.
\end{abstract}

\maketitle

\section{Introduction}
The union-closed sets conjecture is a celebrated open problem in combinatorics which was popularized by Frankl in the late 1970's \cite{Frankl1983}, and is thus often referred to as the Frankl conjecture. Before stating the conjecture, we need a few definitions. Throughout this paper, we think of a \emph{family of sets} or \emph{set system} $\mathcal{F}=(E(\mathcal{F}), \mathcal{S}(\mathcal{F}))$ as being a collection $\mathcal{S}(\mathcal{F})$ of distinct sets $S$ such that every set $S\subseteq E(\mathcal{F})$ where $E(\mathcal{F})$ is the ground set of elements. In general, we let $E(\mathcal{F})=\{1,\ldots, n\}=:[n]$.  A family of sets is said to be \emph{union-closed} if and only if the union of two sets of the family is also a set of the family.

\begin{conjecture}[Frankl, 1979]\label{franklconj}
In a union-closed family $\mathcal{F}$ such that $\mathcal{S}(\mathcal{F})\neq \{\emptyset\}$, there exists an element of $E(\mathcal{F})$ that is in at least half of the sets of $\mathcal{S}(\mathcal{F})$. 
\end{conjecture}

Since 1979, the conjecture has attracted the attention of both lattice theorists as well as combinatorial probabilists, and, more recently, computer scientists.  To the best of our knowledge, this is the first time the problem is investigated through combinatorial optimization. After defining a few more concepts, we give an overview of the literature and present our contributions. 

Let $m(\mathcal{F})$ and $n(\mathcal{F})$ be respectively the numbers of sets and elements of a family $\mathcal{F}$, i.e., $m(\mathcal{F})=|\mathcal{S}(\mathcal{F})|$ and $n(\mathcal{F})=|E(\mathcal{F})|$. Moreover, let $m_e(\mathcal{F})$ be the number of sets in $\mathcal{F}$ containing some element $e\in E(\mathcal{F})$. Let the \emph{degree} of $\mathcal{F}$, denoted by $a(\mathcal{F})$, be the maximum number of sets in $\mathcal{F}$ containing any element of $E(\mathcal{F})$, that is, $a(\mathcal{F})=\max_{e\in E(\mathcal{F})} m_e(\mathcal{F})$. Let $e^*(\mathcal{F})$ be an arbitrary element of maximum degree, i.e., any of possibly many elements in $E(\mathcal{F})$ contained in $a(\mathcal{F})$ sets. For example, if $\mathcal{F}$ is such that $E(\mathcal{F})=[3]$ and $\mathcal{S}(\mathcal{F})=\{\{1,2,3\},\{1,2\},\{1,3\},\emptyset \}$, then $\mathcal{F}$ is union-closed, $m(\mathcal{F})=4$, $n(\mathcal{F})=3$, $m_1(\mathcal{F})=3$, $m_2(\mathcal{F})=2$, $m_3(\mathcal{F})=2$, $a(\mathcal{F})=3$ and $e^*(\mathcal{F})=1$. Since $\frac{3}{4} \geq \frac{1}{2}$, the Frankl conjecture holds for $\mathcal{F}$.

\subsection{A bit of history} Conjecture \ref{franklconj} is known to hold for certain specific families. For example, it has long been known that the conjecture is trivially true for any family $\mathcal{F}$ containing a singleton or a pair, i.e., when there exists $S\in \mathcal{S}(\mathcal{F})$ such that $|S|=1$ or $2$, or when the average size of the sets, $\frac{1}{m(\mathcal{F})}\sum_{S\in \mathcal{S}({\mathcal{F}})} |S|$, is greater or equal to $\frac{n(\mathcal{F})}{2}$. Another early result is the following:

\begin{theorem}[Roberts, 1992]\label{Roberts1992}
The inequality $m(\mathcal{F})<4n(\mathcal{F})-1$ holds for any union-closed family $\mathcal{F}$ that is a minimum counterexample to the Frankl conjecture.
\end{theorem}

The above results and a few others allowed the conjecture to be proven for increasing values of $m$ and $n$ over time (\cite{SarvateRenaud1989}, \cite{SarvateRenaud1990}, \cite{LoFaro19942}, \cite{Poonen1992}, \cite{GaoYu1998}, \cite{Roberts1992}). The current status is as follows:

\begin{theorem}[Roberts and Simpson, 2010]
The Frankl conjecture is true for any family $\mathcal{F}$ with $m(\mathcal{F})\leq 46$.
\end{theorem}

\begin{theorem}[Bo\v{s}njak and Markovi\'c, 2008]
The Frankl conjecture is true for any family $\mathcal{F}$ with $n(\mathcal{F})\leq 11$.
\end{theorem}

More recently, Vu\v{c}kovi\'c and \v{Z}ivkovi\'c announced the following result which is still unpublished.

\begin{theorem}[Vu\v{c}kovi\'c and  \v{Z}ivkovi\'c, 2012]
The Frankl conjecture is true for any family $\mathcal{F}$ such that $n(\mathcal{F})\leq 12$ and $m(\mathcal{F})\leq 50$. 
\end{theorem}

Thus the conjecture is still open for $m(\mathcal{F})\geq 51$ and $n(\mathcal{F})\geq 13$.  A breakthrough in the field is the following result from Reimer \cite{Reimer2003}.

\begin{theorem}\label{ReimerThm}
For any union-closed family $\mathcal{F}$, $\frac{1}{m(\mathcal{F})} \sum_{S\in \mathcal{S}(\mathcal{F})} |S|$, the average size of sets, is at least $\frac{1}{2}\log_2(m(\mathcal{F}))$. 
\end{theorem}

We now turn our attention to three important results which are quite useful for the purposes of this paper. Firstly, Balla, Bollob\'as and Eccles recently proved that the Frankl conjecture holds for families $\mathcal{F}$ containing at least $\frac{2}{3}$ of the sets in the power set of $n(\mathcal{F})$.

\begin{theorem}[Balla, Bollob\'as \& Eccles, 2013]\label{BBE}
The union-closed conjecture holds for any family $\mathcal{F}$ where $m(\mathcal{F}) \geq \frac{2}{3} 2^{n(\mathcal{F})}$.
\end{theorem}

Even more recently, Eccles strengthened this result by proving a stability version in \cite{Eccles2015}. Secondly, instead of proving the Frankl conjecture, one could instead try to prove that any union-closed family contains an element present in at least some fraction of the sets, just as Knill did in the following theorem.

\begin{theorem}[Knill, 1994]\label{knillaverage}
In any union-closed family $\mathcal{F}$, there always exists an element present in at least $\frac{m(\mathcal{F})-1}{\log_2 m(\mathcal{F})}$ sets, that is, $a(\mathcal{F}) \geq \frac{m(\mathcal{F})-1}{\log_2 m(\mathcal{F})}$.
\end{theorem}

W\'ojcik improved this result slightly in \cite{Wojcik1999}, but, amazingly, still no constant fraction is known. Thirdly, Bruhn and Schaudt in \cite{BruhnSchaudt2014} observed the following corollary to Reimer's theorem and the bounds from \cite{VuckovicZivkovic2012}.

\begin{theorem}[Bruhn \& Schaudt, 2013] \label{BSsurvey}
Let $\mathcal{F}$ be any union-closed family $\mathcal{F}$ such that $2^{n(\mathcal{F})-1} < m(\mathcal{F}) \leq 2^{n(\mathcal{F})}$. Then $a(\mathcal{F}) \geq \frac{6}{13}\cdot m(\mathcal{F})$, i.e., there exists an element in a least $\frac{6}{13}$ of the sets of the family.
\end{theorem}

Many other results have been discovered throughout the years. For a more complete history of the problem, we refer the reader to the following excellent survey \cite{BruhnSchaudt2014}. Finally, we note that Timothy Gowers recently led a polymath project, FUNC, on this topic. 

\subsection{Our contributions} In this paper, we examine the Frankl conjecture through a different lens by viewing it as an optimization problem. Indeed, we can rewrite Conjecture \ref{franklconj} either as a maximization or minimization problem, depending on whether we fix $m(\mathcal{F})$ or $a(\mathcal{F})$.   

\begin{conjecture}[Maximization Version]\label{maxconj}
For any positive integer $a$, let $$\mathscr{F}(a)=\{\mathcal{F}|\mathcal{F} \textup{ is a union-closed family}, \mathcal{S}(\mathcal{F})\neq \emptyset \textup{ and } a(\mathcal{F})\leq a\}.$$ Then $\max_{\mathcal{F}\in \mathscr{F}(a)} m(\mathcal{F})\leq 2a$ for all $a\in \mathbb{N}^+$.
\end{conjecture}

\begin{conjecture}[Minimization Version]\label{minconj}
For any positive integer $m$, let $$\mathscr{G}(m)=\{\mathcal{F}|\mathcal{F} \textup{ is a union-closed family}, \mathcal{S}(\mathcal{F}) \neq \emptyset \textup{ and } m(\mathcal{F})=m\}.$$ Then $\min_{\mathcal{F}\in \mathscr{G}(m)} a(\mathcal{F}) \geq \frac{m}{2}$ for all $m\in \mathbb{N}^+$.
\end{conjecture}

Note that the conjectures \ref{franklconj}, \ref{maxconj} and \ref{minconj} are equivalent since there exists a counterexample to the original Frankl conjecture if and only if there exists a union-closed family $\mathcal{F}$ such that $m(\mathcal{F})>0$ and $a(\mathcal{F})<\frac{m(\mathcal{F})}{2}$.

In Section \ref{sec:model}, we model the optimization versions of the Frankl conjecture as integer programs for a fixed $n$ (the number of elements in $E(\mathcal{F})$ for all families $\mathcal{F}$ considered). In Proposition~\ref{fgproperties}, we discuss some of the properties of the optimal values of said programs. Then in Section \ref{sec:compconj}, we present computational results for the models. We observe that the optimal values we computed do not vary as $n$ increases, i.e., 

\begin{align*}
\max_{\substack{\mathcal{F}\in \mathscr{F}(a):\\ n(\mathcal{F})=n}} m(\mathcal{F})= \max_{\substack{\mathcal{F}\in \mathscr{F}(a):\\ n(\mathcal{F})=n+1}} m(\mathcal{F})
\end{align*}
 and 
\begin{align*}
\min_{\substack{\mathcal{F}\in \mathscr{G}(m):\\ n(\mathcal{F})=n}} a(\mathcal{F})=\min_{\substack{\mathcal{F}\in \mathscr{G}(m):\\ n(\mathcal{F})=n+1}} a(\mathcal{F})
\end{align*}
for $n \geq \log_2(a)$. We did not expect this and this is not necessary for the Frankl conjecture to hold. We formally present these two observations as conjectures and prove that they are equivalent. However, these new conjectures do not imply the Frankl conjecture, and conversely, the Frankl conjecture does not imply these two new conjectures. Still, in Section \ref{sec:consequences}, we discuss some of the important implications the new conjectures have on the Frankl conjecture. Notably, proving these conjectures would prove Conjecture \ref{franklconj} for infinitely many values of $|\mathcal{S}(\mathcal{F})|=m$. Moreover, their proof would yield that there always exists an element in $\frac{6}{13}$ of the sets of a union-closed family, and would thus achieve the first known constant bound for the percentage of sets containing some element in any union-closed family. Finally, in Section \ref{sec:newresults}, we prove a restricted version of the new conjectures using an observation of Falgas-Ravry in \cite{FalgasRavry2011} and discuss the importance of twin sets, which we define as sets that differ only in one element.

\section{The Frankl Integer Problems and Two New Conjectures}\label{sec2}

\subsection{Modeling the Frankl optimization problems}\label{sec:model}

For any positive integers $a, n$, let $$\mathscr{F}(n,a)=\{\mathcal{F}\in \mathscr{F}(a)|n(\mathcal{F})=n\}.$$ Then proving that $\max_{\mathcal{F}\in \mathscr{F}(n,a)} m(\mathcal{F}) \leq 2a$ for all possible $n, a$ would prove Conjecture \ref{maxconj} (and thus the original conjecture). Fix $n, a$, and let $f(n,a)=\max_{\mathcal{F} \in \mathscr{F}(n,a)} m(\mathcal{F})$. Then we can find $f(n,a)$ by solving the following integer program

\boxedeq{%
f(n,a)= \max & \sum_{S\in \mathcal{S}_n} x_S &&  \\
\textup{such that } & x_U+x_T \leq 1+x_S \quad \quad\quad \quad\quad \quad \hspace{0.05cm} \forall T\cup U = S \in \mathcal{S}_n \\ 
& \sum_{S\in \mathcal{S}_n: e\in S} x_S \leq a \quad \quad\quad \quad\quad \quad \quad \forall e\in [n] \\ 
& x_S \in \{0,1\}\quad \quad\quad \quad\quad\quad\quad\quad \quad \hspace{0.1cm}\forall S\in \mathcal{S}_n,
}
where $\mathcal{S}_n$ is the power set of $[n]$, and the variable $x_S$ for any set $S\in \mathcal{S}_n$ is $1$ if $S$ is in the family, and $0$ otherwise. Thus, we maximize the number of sets while ensuring the family is union-closed (through the first constraint) and that $a(\mathcal{F}) \leq a$ holds (through the second constraint), i.e., we calculate $f(n,a)$. 

Similarly, let $$\mathscr{G}(n,m)=\{\mathcal{F}\in \mathscr{G}(m)|n(\mathcal{F})=n\}$$ for any positive integers $m$ and large enough $n$. Then, it is easy to see that proving that $\min_{\mathcal{F} \in \mathscr{G}(n,m)} a(\mathcal{F}) \geq \frac{m}{2}$ for all large enough $n$ and $m$ would prove Conjecture \ref{minconj}. Fix $n, m$, and let $g(n,m)=\min_{\mathcal{F}\in \mathscr{G}(n,m)} a(\mathcal{F})$ for $m \geq 2$. Then we can find $g(n,m)$ by solving the following integer program

\boxedeq{%
g(n,m) = \min & \sum_{S\in \mathcal{S}_n: 1\in S} x_S  \\
\textup{such that } &  x_U+x_T \leq 1+x_S \hspace{0.2cm} \quad \quad \quad \quad \quad \forall T\cup U = S \in \mathcal{S}_n \\
& \sum_{S\in \mathcal{S}_n: i\in S} x_S \geq \sum_{S\in \mathcal{S}_n: j\in S} x_S \quad \quad \forall 1\leq i\leq j\leq n\\
& \sum_{S\in \mathcal{S}_n} x_S = m  \\
& x_S \in \{0,1\} \hspace{0.25cm}\quad\quad\quad\quad\quad\quad\quad\quad \forall S\in \mathcal{S}_n,
}
where the variables $x_S$ are as before. The second constraint says that element $1$ is the element contained in the most number of sets in the family, so we are minimizing the maximum number of sets containing the most frequent element while enforcing that the family is union-closed (through the first constraint) and has $m$ sets (through the third constraint). 

Most work done on the Frankl conjecture has been from the $g(n,m)$ point of view, not $f(n,a)$. Moreover, in \cite{Renaud1991} and \cite{Renaud1995}, Renaud defined $\varphi(m)$ to be the minimum number of sets containing the most frequent element in a family among all union-closed families on $m$ sets. Some of the properties he proves for $\varphi(m)$ are not unlike those we will prove for $g(n,m)$.

Note moreover that there are values of $a, m, n$ for which $f(n,a)$ and $g(n,m)$ have trivial solutions that do not interest us. For example, it is clear that $f(n,a)=2^n$ if $a\geq 2^{n-1}$. Indeed, the power set of $n$, $\mathcal{S}_n$, is union-closed, and there are $2^n$ sets in $\mathcal{S}_n$ where each element is in exactly $2^{n-1}$ sets. It is thus a trivially optimal solution for $f(n,a)$. It is also clear that $g(n,m)$ has trivially no solution if $m>2^n$. Indeed, even if we take all of the $2^n$ sets in $\mathcal{S}_n$, we would have less sets than the number of sets required by the program. We first study a few properties of the functions $f$ and $g$ for non-trivial values of $a, m, n$. 

\begin{proposition}\label{fgproperties} 
The following properties hold.
\begin{enumerate}
\item The function $f$ is non-decreasing in $n$, that is, $f(n,a) \leq f(n+1,a)$ for every $a, n \in \mathbb{N}^+$ such that $n \geq \lceil \log_2 a \rceil + 1$.\label{fnondecreasing}
\item The function $g$ is non-increasing in $n$, that is, $g(n,m) \geq g(n+1,m)$ for every $m, n\in \mathbb{N}^+$ such that $n \geq \lceil \log_2 m \rceil$.\label{gnonincreasing}
\item The function $f$ is strictly increasing in $a$, that is, $f(n,a) < f(n,a+1)$ for every $a, n\in \mathbb{N}^+$ such that $n > \lceil \log_2 a \rceil + 1$.\label{fstrict}
\item The function $g$ is non-decreasing in $m$, that is, $g(n,m) \leq g(n,m+1)$ for every $m,n\in \mathbb{N}^+$ such that $n \geq \lceil \log_2 m \rceil$.\label{gnondecreasing}
\item We have that $g(n,f(n,a)) = a$ for all $a, n\in \mathbb{N}^+$ such that $n > \lceil \log_2 a \rceil + 1$.\label{gfequal}
\item We have that $f(n,g(n,m)) \geq m$ for all $m,n \in \mathbb{N^+}$ such that $n \geq \lceil \log_2 m \rceil$.\label{fggeq}
\end{enumerate}
\end{proposition}

\proof
\begin{enumerate}
\item Fix $a, n$ and suppose $f(n,a) = m$. Take a family $\mathcal{F}$ that is optimal, i.e., a family $\mathcal{F}\in \mathscr{F}(n,a)$ such that $m(\mathcal{F})=m$. Add an element $n+1$ to the family such that this element is exactly in the same sets as some other element $e$ of the family. Then this augmented family is still union-closed, and every element of it is still in at most $a$ sets. Thus, $f(n+1,a) \geq m = f(n,a)$.

\item Just as for (\ref{fnondecreasing}), we clone an element. 

\item Fix $a, n$ and suppose $f(n,a)=m$. Take a family $\mathcal{F}$ that is optimal (as before). Then add to this family one of the largest sets that is not already present in the family, i.e., a set in $\mathcal{S}_n\backslash \mathcal{S}(\mathcal{F})$ containing as many elements as possible. This is always possible if $f(n,a)<2^n$ (note that this is the case if $a<2^{n-1}$ since we cannot take all of $\mathcal{S}_n$ then). The new family we built has $m+1$ sets and is still union-closed since taking the union of any other set with the added set will give either the new set itself or a greater set (which is present in the family by construction). Moreover, every element in this new family is present in at most $a+1$ sets. Thus, $f(n,a+1)\geq m+1 = f(n,a) + 1$.

\item Fix $m,n$ and suppose $g(n,m+1)=a$. Take a family $\mathcal{F}$ that is optimal (as before). Remove the smallest set of the family, i.e., one of the sets in $\mathcal{S}(\mathcal{F})$ with the least number of elements. The family is still union-closed since the set that was removed was not the union of any other two sets since they would have to be smaller. Moreover, this new family has $m$ sets and every element is still there at most $a$ times, so $g(n,m) \leq a=g(n,m+1)$.

\item Fix $a, n$ and suppose that $f(n,a)=m$. Thus $m$ is the maximum number of sets in any union-closed family $\mathcal{F}$ on $n$ elements with $a(\mathcal{F})\leq a$. Because $m \leq 2^n$, $g(n,m)$ is feasible. Certainly, this implies that $g(n,m)=g(n,f(n,a)) \leq a$. Suppose that $g(n,f(n,a))=a' < a$. Then $f(n,a') \geq m$. Since we have already shown that the function $f$ strictly increases in $a$ when $n>\lceil \log_2 a \rceil +1$, this is a contradiction.

\item Suppose that $g(n,m)=a$. Then $n \geq \lceil \log_2 m \rceil$, else $g(n,m)$ would be infeasible. Thus, the most frequent element in a union-closed family with $n$ elements and $m$ sets is present in a least $a$ sets. Certainly, this means that $f(n,a)=f(n,g(n,m))\geq m$.
\end{enumerate}
\qed

\subsection{Computations and Conjectures}\label{sec:compconj}

We computed $f(n,a)$ and $g(n,m)$ for different values with the mixed-integer commercial solver IBM ILOG CPLEX version 12.4. See http://www.math.washington.edu/$\sim$raymonda/frankl.py for the source code to generate the .lp files. Table \ref{tablefg} contains some of the results we obtained.
\begin{table}[h!]
 \caption{Values of $f(n,a)$ and $g(n,m)$ (respectively left and right)}
 \label{tablefg}
\begin{center}
\setlength{\tabcolsep}{2.6pt}
\begin{tabular}{|l|cccccccc|}
\hline
$a\backslash n$ & 1 & 2 & 3 & 4 & 5& 6& 7& 8\\
  \hline
  1 & 2 & 2 & 2 & 2 & 2& 2& 2& 2\\
  2 & \textcolor{gray}{2} & 4 & 4 & 4& 4& 4& 4& 4\\
  3 & \textcolor{gray}{2}& \textcolor{gray}{4}& 5 & 5 & 5 & 5 & 5 & 5\\
  4 & \textcolor{gray}{2}& \textcolor{gray}{4}& 8& 8& 8& 8& 8& 8 \\ 
  5 & \textcolor{gray}{2}& \textcolor{gray}{4}& \textcolor{gray}{8}& 9& 9& 9& 9& 9\\
  6 & \textcolor{gray}{2}& \textcolor{gray}{4}& \textcolor{gray}{8}& 10& 10& 10& 10& 10 \\
  7 & \textcolor{gray}{2}& \textcolor{gray}{4}& \textcolor{gray}{8}& 12& 12& 12& 12& 12 \\
  8 & \textcolor{gray}{2}& \textcolor{gray}{4}& \textcolor{gray}{8}& 16& 16& 16& 16& 16 \\
  9 & \textcolor{gray}{2}& \textcolor{gray}{4}& \textcolor{gray}{8}& \textcolor{gray}{16}& 17& 17& 17& 17 \\
  10 & \textcolor{gray}{2}& \textcolor{gray}{4}& \textcolor{gray}{8}& \textcolor{gray}{16}& 18& 18& 18& 18 \\
  11 & \textcolor{gray}{2}& \textcolor{gray}{4}& \textcolor{gray}{8}& \textcolor{gray}{16}& 19& 19& 19& 19  \\
  12 & \textcolor{gray}{2}& \textcolor{gray}{4}& \textcolor{gray}{8}& \textcolor{gray}{16}& 21& 21& 21& 21 \\
  13 & \textcolor{gray}{2}& \textcolor{gray}{4}& \textcolor{gray}{8}& \textcolor{gray}{16}& 23& 23& 23& 23 \\
  14 & \textcolor{gray}{2}& \textcolor{gray}{4}& \textcolor{gray}{8}& \textcolor{gray}{16}& 25 & 25 & 25 & 25 \\
  15 & \textcolor{gray}{2}& \textcolor{gray}{4}& \textcolor{gray}{8}& \textcolor{gray}{16}& 27& 27& 27& 27 \\
  16 & \textcolor{gray}{2}& \textcolor{gray}{4}& \textcolor{gray}{8}& \textcolor{gray}{16}& 32& 32& 32& 32 \\
  \hline
\end{tabular} \quad \begin{tabular}{|l|cccccccc|}
\hline
$m\backslash n$  & 1 & 2 & 3 & 4 & 5& 6& 7& 8 \\
  \hline
   &  &  &  &  & & & &   \\
  2 & 1& 1 &1 & 1& 1& 1& 1& 1 \\
  3 &\textcolor{gray}{-} & 2& 2 & 2 & 2 & 2 & 2 & 2 \\
  4 &\textcolor{gray}{-} & 2& 2& 2& 2& 2& 2& 2\\ 
  5 &\textcolor{gray}{-} &\textcolor{gray}{-} &3 & 3& 3& 3& 3& 3\\
  6 &\textcolor{gray}{-} &\textcolor{gray}{-} & 4& 4& 4& 4& 4& 4 \\
  7 &\textcolor{gray}{-} &\textcolor{gray}{-} & 4& 4&4& 4& 4& 4\\
  8 &\textcolor{gray}{-} &\textcolor{gray}{-} & 4& 4& 4& 4& 4& 4\\
  9 &\textcolor{gray}{-} &\textcolor{gray}{-} &\textcolor{gray}{-} & 5& 5& 5& 5& 5\\
  10 &\textcolor{gray}{-} &\textcolor{gray}{-} &\textcolor{gray}{-} & 6& 6& 6& 6& 6\\
  11 &\textcolor{gray}{-} &\textcolor{gray}{-} & \textcolor{gray}{-}& 7& 7& 7& 7& 7\\
  12 &\textcolor{gray}{-} &\textcolor{gray}{-} & \textcolor{gray}{-}& 7& 7& 7& 7& 7\\
  13 &\textcolor{gray}{-} &\textcolor{gray}{-} & \textcolor{gray}{-}& 8& 8& 8& 8& 8\\
  14 &\textcolor{gray}{-} &\textcolor{gray}{-} & \textcolor{gray}{-}& 8& 8 & 8 & 8 & 8\\
  15 &\textcolor{gray}{-} &\textcolor{gray}{-} & \textcolor{gray}{-}& 8& 8& 8& 8& 8\\
  16 &\textcolor{gray}{-} &\textcolor{gray}{-} & \textcolor{gray}{-}& 8& 8& 8& 8& 8\\
\hline
\end{tabular}
\end{center}
\end{table}

 A clear pattern emerges: For any $n\geq\lceil \log_2 a \rceil + 1$, that is, for any non-trivial value of $n$ when $a$ is fixed, $f(n,a)$ takes the same value as $n$ increases.  Similarly, for any  $n \geq \lceil \log_2 m \rceil$, that is, for any non-trivial value of $n$ when $m$ is fixed, $g(n,m)$ takes the same value as $n$ increases.

To the best of our knowledge, this has never been observed before. We formulate these observations as conjectures.

\begin{conjecture}[$f$-conjecture]\label{fconjecture}
Fix $a\in \mathbb{N}^+$.  Then $f(n,a) = f(n+1,a)$ for every $n\in \mathbb{N}^+$ such that $n \geq \lceil \log_2 a \rceil + 1$.
\end{conjecture}

\begin{conjecture}[$g$-conjecture]\label{gconjecture}
Fix $m\in \mathbb{N}^+$. Then $g(n,m) = g(n+1,m)$ for every $n\in \mathbb{N}^+$ such that $n \geq \lceil \log_2 m  \rceil $.
\end{conjecture}

We checked these conjectures computationally up to $n=9$ for all non-trivial values of $a$ for $f(n,a)$ and up to $n=8$ for all non-trivial values of $m$ for $g(n,m)$ (see Appendix).  We first show that these two conjectures are equivalent.

\begin{theorem}\label{conjecturesequivalent}
We have that $f(n,a)=f(n+1,a)$ for every $a,n\in \mathbb{N}^+$ such that $n \geq \lceil \log_2 a \rceil + 1$ if and only if $g(n',m)=g(n'+1,m)$ for every $m,n' \in \mathbb{N}^+$, $m \geq 2$ such that $n' \geq \lceil \log_2 m \rceil$.
\end{theorem}

\proof
Suppose that $g(n',m)=g(n'+1,m)$ for every $m,n' \in \mathbb{N}^+$, $m \geq 2$ such that $n' \geq \lceil \log_2 m \rceil$.  Pick an arbitrary $a\in \mathbb{N}^+$ and choose any $n\in \mathbb{N}^+$ such that $n\geq \lceil \log_2 a \rceil + 1$. By Proposition~\ref{fgproperties}(\ref{gfequal}), we know that there exists $m$ such that $g(n,m)=a$, namely $m:=f(n,a)$. Let $m^*$ be the greatest number for which $g(n,m^*)=a$. By definition, it follows that $f(n,a) \leq m^*$. Moreover,  $f(n,a) \geq m^*$ since $f(n,g(n,m^*)) \geq m^*$ by Proposition~\ref{fgproperties}(\ref{fggeq}). Thus, $f(n,a) = m^*$ for any $n\geq \lceil \log_2 a \rceil + 1$.

Suppose that $f(n,a)=f(n+1,a)$ for every $a,n\in \mathbb{N}^+$ such that $n \geq \lceil \log_2 a \rceil + 1$. Pick an arbitrary $m\geq 2$, and choose any $n'\geq \lceil \log_2 m \rceil$. If there exists $a$ such that $f(n',a)=m$, then, by Proposition~\ref{fgproperties}(\ref{gfequal}), $g(n',m)=g(n',f(n',a))=a$. Therefore, $g(n',m)=g(n'+1,m)$ for all $n' \geq \lceil \log_2 m \rceil$. If there does not exists $a$ such that $f(n',a)=m$, then, by Proposition~\ref{fgproperties}(\ref{fggeq}), $f(n',g(n',m))=m+b$, $b>0$, and $g(n',m)=g(n',m+b)$. Then $g(n', m+b)=:a'$ for all $n'$ since $f(n'(g(n', m+b))=m+b$, which brings us back to the first case. Therefore, $g(n',m)=g(n',m'+b)=a'$ for all $n'\geq \lceil \log_2 m \rceil$. 
\qed

The two conjectures are thus equivalent: proving one would prove the other. Another way to view the $f$-conjecture is as follows: to construct an optimal solution for $f(n,a)$, first construct an optimal solution for $f(\lceil \log_2 a \rceil + 1, a)$, and make $n-\lceil \log_2 a \rceil -1$ copies of some element, i.e., put these new elements in the same sets as the original element, or if you prefer, put these new elements in none of the sets. If the $f$-conjecture is true, such families would be optimal for $f(n,a)$. Note though that there can also be other optimal families: the conjecture simply states that families obtained through this process are optimal. The same idea can be applied to the $g$-conjecture. 

Note that these conjectures are different from the Frankl conjecture. For one thing, even if the $f$- and $g$-conjectures hold, one would still need to show that $f(\lceil \log_2 a \rceil + 1,a) \leq 2a$ for every $a$ or that $g(\lceil \log_2 m  \rceil,m) \geq \frac{m}{2}$ for all $m$ to prove the Frankl conjecture; therefore the $f$- and $g$-conjectures do not imply the Frankl conjecture. Moreover, the Frankl conjecture does not immediately imply the $f$- and $g$-conjectures.  Certainly, if the Frankl conjecture is true, then 
\begin{align*}
\max_{n} f(n,a) \leq 2a
\end{align*}
for all $a\in \mathbb{N}^+$, else the Frankl conjecture would not be true; however, how the function $f(n,a)$ behaves for different $n$ for a fixed $a$ is irrelevant. Certainly, one can easily observe that the Frankl conjecture implies that the function $f$ has to stabilize at some point as $n$ increases since, by Proposition~\ref{fgproperties}(\ref{fnondecreasing}), $f(n,a)\leq f(n+1,a)$, and since  $f(n,a)\in \mathbb{N}$ by definition. However, the Frankl conjecture does not imply that the function $f$ should be stable immediately as $n\geq \lceil \log_2 a \rceil + 1$, that is, as soon as there are enough elements for there to be at least $a$ sets containing an element $e$.  Similarly, if the Frankl conjecture is true, then
\begin{align*}
\min_{n\geq \lceil \log_2 m \rceil} g(n,m) \geq\frac{m}{2}
\end{align*}
for all $m\in \mathbb{N}^+$, and so again one can observe that the $g$ function has to stabilize at some point as $n$ increases since we show in Proposition~\ref{fgproperties}(\ref{gnonincreasing}) that $g(n,m) \geq g(n+1,m)$. But again, the Frankl conjecture does not imply that the function $g$ has to stabilize immediately when $n \geq \lceil \log_2 m\rceil$, that is, as soon as there are enough elements for there to be at least $m$ sets in the power set of $n$.

Therefore, the $f$- and $g$-conjectures do not imply the Frankl conjecture and the latter does not imply the former two either. Still, proving the $f$- and $g$-conjectures would have wide-reaching implications for the Frankl conjecture.

\subsection{Consequences of the $f$- and $g$-conjectures on the Frankl conjecture}\label{sec:consequences}

Note that if the $f$- and $g$-conjectures are true, then proving the union-closed sets conjecture for \emph{large} families, i.e., families such that $m\geq 2^{n-1}+1$ would be enough to prove it for \emph{all} families. Such a proof does not yet exist, however, combining the $f$- and $g$-conjectures with Theorem \ref{BBE} would already go a long way towards solving the Frankl conjecture. 

\begin{theorem}\label{fgBBE}
If the $f$- and $g$-conjectures hold, then Conjecture \ref{minconj} holds for all $m$ for which there exists $i\in \mathbb{N}^+$ such that $\frac{2}{3}2^i \leq m \leq 2^i$.
\end{theorem}

\proof Let $m$ be such that $\frac{2}{3}2^i \leq m \leq 2^i$ for some $i$. If $g(n,m)<\frac{m}{2}$ for some $n\geq \lceil \log_2 m \rceil$, then we have a counterexample: a union-closed family on $n$ elements and $m$ sets where every element is in less than half of the sets. By the $g$-conjecture, since $i \geq \lceil \log_2 m \rceil$, $g(n,m)=g(i, m)<\frac{m}{2}$. However, by Theorem \ref{BBE}, $g(i,m) \geq \frac{m}{2}$ since $m\geq \frac{2}{3}2^i$. Thus, we have obtained a contradiction and $g(n,m)\geq \frac{m}{2}$ for all $n\geq \lceil \log_2 m \rceil$. 
\qed

This would mean that Conjecture \ref{minconj} would hold for about $\frac{2}{3}$ of all possible values of $m$. Recall that Conjecture \ref{minconj} is equivalent to the Frankl conjecture, and so Theorem~\ref{fgBBE} could be reformulated in such terms. 

Another nice consequence of the $f$- and $g$-conjectures would be that there would finally be a known constant fraction of sets containing $e^*(\mathcal{F})$, the most frequent element in $\mathcal{F}$, in any union-closed family $\mathcal{F}$.

\begin{theorem}
If the $f$- and $g$-conjectures hold, then any union-closed family on $m$ sets contains an element in at least $\frac{6}{13} m$ sets of the family.
\end{theorem}

\proof
By Theorem \ref{BSsurvey}, we know that $g(\lceil \log_2 m  \rceil,m) \geq \frac{6}{13}m$ for any $m\in \mathbb{N}^+$. By the $g$-conjecture, we know that $g(n,m) \geq \frac{6}{13}m$ for all $n\geq \lceil \log_2 m  \rceil$. Therefore, we know that any family on $m$ sets contains an element in at least  $\frac{6}{13}m$ sets of the family.
\qed

Thus, from our point of view, studying the $f$- and $g$-conjectures offers new ways of attacking the Frankl conjecture, in addition to being interesting in and of itself. Hence, the new conjectures warrant a closer examination.

\section{Towards Proving the New Conjectures}\label{sec:newresults}
\subsection{Twin sets and a partial proof of the new conjectures}
As noted before in Theorem \ref{conjecturesequivalent}, the $f$- and $g$-conjectures are equivalent, so from now on, we will focus only on the $f$-conjecture. We now introduce a new idea: twin sets.

\begin{definition}
We call two sets $S_1,S_2$ with $n>|S_1|>|S_2|$ \emph{twin sets} if $|S_1\triangle S_2| = 1$ . We call $S_1$ the \emph{big twin} and $S_2$ the \emph{little twin}. Moreover, we call the element $e=S_1\triangle S_2$ the \emph{twin difference of $S_1$ and $S_2$}.
\end{definition}

Twin sets play an important role in proving the $f$-conjecture as made clear by the following lemma.

\begin{lemma}\label{ftwinscont}
Suppose that the $f$-conjecture is not true and that there exists values of $n$ and $a$ with $n \geq \lceil \log_2 a \rceil + 1$ such that $f(n,a)<f(n+1,a)$. Then every optimal solution of $f(n+1,a)$ is such that every element is a twin difference for at least one pair of sets. 
\end{lemma}

\proof
Suppose that for $f(n+1,a)=:m$ there exists an optimal solution for which there exists an element that is not a twin difference. Then removing this element (from every set containing it) will leave all of the sets distinct, and the family will still be union-closed, but with $n$ elements. In that case, we know that $f(n,a) \geq m = f(n+1,a)$, and since $f(n,a) \leq f(n+1,a)$ by Proposition~\ref{fgproperties}(\ref{fnondecreasing}), we obtain $f(n,a) = f(n+1,a)$ as in the conjecture, a contradiction.
\qed

Therefore we only need to focus on the case where every optimal solution for some $f(n,a)$ is such that every element is a twin difference.\\

Using Lemma \ref{ftwinscont} and an observation from Falgas-Ravry \cite{FalgasRavry2011}, we can prove the $f$-conjecture for the cases when $n>a$. 

\begin{theorem}\label{falgasravryf}
We have that $f(n-1,a)=f(n,a)$ for all $n> a$.
\end{theorem}

\proof
First assume $f(n-1,a)<f(n,a)$ for some $n> a$. By Lemma \ref{ftwinscont}, every element in any optimal solution of $f(n,a)$ is a twin difference of at least one pair of sets. This implies that any optimal solution for $n$ elements contains no two elements that are exactly in the same sets. Indeed, such elements would not be twin differences. We can thus apply the same construction as in \cite{FalgasRavry2011}. Order the elements $[n]$ by decreasing frequency. Observe now that for all elements $1\leq i<j\leq n$, there exists $S_{ij}$ such that $i\in S_{ij}$ and $j\not\in S_{ij}$. For $2\leq j\leq n$, let $S_j=\cup_{i=1}^{j-1} S_{ij}$, and let $S_{n+1}=[n]$. Note that the $S_j$'s are all distinct since $[1,j-1] \subseteq S_j$ and $j\not \in S_j$. Since the most frequent element, element $1$, is in at least these $n$ sets of the family, we have that $a\geq n$, which is a contradiction.
\qed

First note that the Falgas-Ravry construction can also be applied to the function $g$ to prove that $g(n,m)=g(n+1,m)$ for all $n\geq m-1$. Furthermore, observe that having no two elements in exactly the same sets is a much weaker constraint than having each element being a twin difference. Therefore, it might be possible to improve this result.

\subsection{Models with twins and some computations}

Since we are interested in union-closed families in which every element is a twin difference, we can modify the models from section \ref{sec:model} so that they only consider such families.

Additionally, we will not count the set with every element as a possible big twin (that we will call trivial twin) as this case also yields that $f(n+1,a)=f(n,a)$. Indeed, if an element $e$ is the twin difference only of the trivial twin and $[n]\backslash \{e\}$, then we can remove $e$ and replace the set $[n]\backslash\{e\}$ (now a set with all $n-1$ elements identical to the trivial twin) with the biggest set missing from the family (for example, the first set missing in the lexicographic order using any ordering of the elements). This new family is union-closed by the same argument as in Proposition~\ref{fgproperties}(\ref{fstrict}), and each element is still in at most $a$ sets, and so here again we have that $f(n+1,a)=f(n,a)$. 

Enforcing that each element is a non-trivial twin difference for some pair of sets is easily done in our $f$- and $g$-programs by introducing the variable $z_S^e$ for every $S\in \mathcal{S}_n$ and $e\in [n]$, which is zero if at least one of $S$ and $S\cup e$ is absent from the family. This of course makes the program much larger. Luckily, $z_S^e$ can be a continuous non-negative variable. Indeed, by adding the constraints $z_S^e \leq x_S$ and $z_S^e \leq x_{S\cup e}$, we ensure that the twin variable is zero if $S$ or $S\cup e$ is missing.  We also add the constraint $\sum_{\substack{S\not\ni e,\\ |S|\neq n-1}} z_S^e \geq 1$ for every $e\in [n]$, which ensures each element is a non-trivial twin difference. We let $f_t(n,a)$ and $g_t(n,m)$ be the optimal values of the programs for $f(n,a)$ and $g(n,m)$ with the new variables and constraints. 

 We present in Table \ref{fnawithtwins} computational results for $f_t(n,a)$ and $g_t(n,m)$.

\begin{table}[h!]
 \begin{center}
\caption{Non-trivial values of $f_t(n,a)$ and $g_t(n,m)$ (respectively left and right)}
\label{fnawithtwins}
\setlength{\tabcolsep}{2.6pt}
\begin{tabular}{|l|cccccccc|}
\hline
 $a\backslash n$ & 1 & 2 & 3 & 4 & 5& 6& 7& 8 \\
  \hline
  1 & $-$ & $-$ & $-$& $-$ & $-$& $-$& $-$& $-$  \\
  2 &  & $-$ & $-$ & $-$& $-$& $-$& $-$& $-$ \\
  3 & & & $-$ & $-$ & $-$ & $-$ & $-$ & $-$ \\
  4 & & & 8& $-$& $-$& $-$& $-$& $-$ \\
  5 & & & & 8& $-$& $-$& $-$& $-$\\
  6 & & & & 10&  9&  $-$&  $-$&  $-$ \\
  7 & & & & 12& 11& 10&  $-$&  $-$\\
  8 & & & & 16& 13& 12& 11& $-$\\
  9 & & & & & 15& 14& 13& 12\\
  10 & & & & & 18& 16& 15& 14\\
  11 & & & & & 19& 19& 17& 16\\
  12 & & & & & 21& 20& 20& 18\\
  13 & & & & & 23& 22& 21& 21\\
  14 & & & & & 25& 24& 23& 22\\
  15 & & & & & 27& 25& 25& 24\\
  16 & & & & & 32& 28& 26& 26\\
  \hline
\end{tabular} \quad \begin{tabular}{|l|ccccccc|}
\hline
 $m\backslash n$ & 1 & 2 & 3 & 4 & 5& 6& 7\\
  \hline
  1 & $-$ & $-$ & $-$ & $-$ & $-$ & $-$ & $-$ \\
  2 & $-$ & $-$ & $-$ & $-$ & $-$ & $-$ & $-$ \\
  3 &     & $-$ & $-$ & $-$ & $-$ & $-$ & $-$ \\
  4 &     & 2   & $-$ & $-$ & $-$ & $-$ & $-$ \\
  5 &     &     & 4   & $-$ & $-$ & $-$ & $-$ \\
  6 &     &     & 4   & 5   & $-$ & $-$ & $-$ \\
  7 &     &     & 4   & 5   & 6   & $-$ & $-$ \\
  8 &     &     & 4   & 5   & 6   & 7   & $-$ \\
  9 &     &     &     & 6   & 6   & 7   & 8   \\
  10 &    &     &     & 6   & 7   & 7   & 8   \\
  11 &    &     &     & 7   & 7   & 8   & 8   \\
  12 &    &     &     & 7   & 8   & 8   & 9   \\
  13 &    &     &     & 8   & 8   & 9   & 9   \\
  14 &    &     &     & 8   & 9   & 9   & 10  \\
  15 &    &     &     & 8   & 9   & 10  & 10  \\
  16 &    &     &     & 8   & 10  & 10  & 11  \\
  \hline
\end{tabular}
\end{center}
\end{table}

From Theorem \ref{falgasravryf}, we know that everything past the main diagonal for $f_t(n,a)$ is infeasible. Similarly, everything past the lower diagonal for $g_t(n,m)$ is infeasible.

The fact that $f_t(n,a)$ decreases as $n$ increases appears counterintuitive at first, but on second thought it makes sense. If there are more elements, each forced to be a twin difference, then there will be more distinct unions of sets. Therefore the risk of violating the $a$-limit increases, and so the number of allowable sets decreases. If one could prove that, then the $f$-conjecture would be proven, and thus the $g$-conjecture as well. Another idea pointing in a similar direction is the following.

\begin{lemma}
For a fixed $a$, the minimum number of twin pairs for any element in an optimal union-closed family for $f(n,a)$ is bounded above by $2(a-n+1)$. In particular, this upper bound decreases as $n$ increases.
\end{lemma}

\proof
Order the elements by decreasing frequency. Let $t_{\mathcal{F}}$ be the minimum number of twin pairs for any element in an optimal family $\mathcal{F}$ for $f(n,a)$. If $t_{\mathcal{F}}\geq 1$, then, by the Falgas-Ravry observation presented in Theorem \ref{falgasravryf}, there exists a set $S_j$ with $[1,j-1]\subseteq S_j$ and $j\not\in S_j$ for every $2\leq j \leq n$, as well as the set $S_{n+1}=[n]$. Note that, for element $1$, these sets can only be big twins (since $1$ is in all of these sets). Suppose $u$ of these sets are big twins for element $1$. Then there exists at least $t_{\mathcal{F}}-u$ other sets that are big twins for element $1$. Thus element $1$ is in at least $n+t_{\mathcal{F}}-u$ sets. Moreover, element $2$ will be in at least $u-1$ of the small twins for element $1$, and it is in $n-1$ sets among the $S_j$'s (which are distinct from these small twins). Therefore, element $2$ is in at least $n-1+u-1$ sets. Therefore, $a\geq \max\{n+t_{\mathcal{F}}-u, n+u-2\} \geq n-1+\frac{t_{\mathcal{F}}}{2}$. This implies that $t_{\mathcal{F}} \leq 2(a-n+1)$. Thus, for a fixed $a$, as $n$ increases, the potential number of twins in an optimal family decreases. 
\qed

Some of the properties in Proposition~\ref{fgproperties} of the $f$- and $g$-functions also hold for $f_t$ and $g_t$.

\begin{proposition}
The following properties hold.
\begin{enumerate}
\item We have that $f_t(n,a) < f_t(n,a+1)$ if $a < 2^{n-1}$.\label{strictt}
\item We have that $g_t(n,f_t(n,a)) = a$ for all $a$ and $n$.
\item We have that $f_t(n,g_t(n,m)) \geq m$ for all $m$ and $n$.
\end{enumerate}
\end{proposition}

\proof
\begin{enumerate}
\item Suppose $f_t(n,a)=m$. Take a family $\mathcal{F}$ that is optimal. Then add to this family one of the greatest set, i.e., a set containing as many elements as possible or lexicographically greatest for some ordering of the elements,  that is not already present in the family. This is always possible since $f_t(n,a)<2^n$. This new family is still union-closed, and moreover there are still twins for each element since we did not remove any set. Each element is now in at most $a+1$ sets, so this family is valid for $f_t(n,a+1)$, so $f_t(n,a+1) \geq f_t(n,a)+1$.
\item Suppose that $f_t(n,a)=m$. This means the maximum number of sets in a union-closed family such that every element is in at most $a$ sets and such that each element is a non-trivial twin difference is $m$. Certainly, this means that $g_t(n,m)=g_t(n,f_t(n,a)) \leq a$ since we're minimizing and the previous family is valid here. Suppose now that $g_t(n,f_t(n,a))=a' < a$. Then $f_t(n,a') \geq m$. Since $f_t$ is strictly increasing in $a$ by (\ref{strictt}), this is a contradiction.
\item Suppose that $g_t(n,m)=a$. This means that there exists a family on $m$ sets where each element is a twin difference and is also present in at most $a$ sets. Thus, this family is valid for $f_t(n,a)$, and so $f_t(n,a) = f_t(n,g_t(n,m)) \geq m$ since we are maximizing.
\end{enumerate}
\qed

\subsection{Number of twin pairs}

Another direction worth investigating would be to prove that if $f(n,a)<f(n+1,a)$, then every element is the difference of an increasingly large number of twins. At some point, this ceases to be possible (trivially, an element cannot be the twin difference of more than $a$ twin pairs), and so we would reach a contradiction.

\begin{theorem}\label{twotwins}
If $f(n,a)<f(n+1,a)$, then every element in a $f(n+1,a)$-optimal family is the difference of at least two pairs of twin sets.
\end{theorem}

\proof
Suppose that $f(n,a)=m$ and $f(n+1,a)=m+k$ for some $k>0$. Then we know that any optimal solution for $f(n+1,a)$ must be such that every element is a twin difference, otherwise we could remove that element and get an $m+k$ union-closed family spanning $n$ elements such that none is in more than $a$ sets, and so $f(n,a)\geq m+k$, a contradiction. Now let

$$k':=\min_{e\in [n+1]}|\{S\in \mathcal{F}|e\not\in S, \  e\cup S\in \mathcal{F} \textup{ and } \mathcal{F} \textup{ is an optimal family for } f(n+1,a)\}|,$$
i.e., $k'$ is the minimum number of twin pairs for which an element is a twin difference in an optimal solution for $f_{n+1}(a)$.

We now show that $k\leq k'$. Suppose not. Let $e'$ and $\mathcal{F}'$ be an element and a family such that $e'$ is a twin difference for $k'$ twin pairs. Remove $e'$ from $\mathcal{F}'$, and remove the $k'$ sets that are now duplicated. Call this new family $\mathcal{F}''$. What remains is a union-closed family of $m+k-k'$ sets on $n$ elements. So $m=f(n,a) \geq m+k-k'$, which implies that $k\leq k'$.

Now suppose that $k=k'$. Notice then that there must exist $e''$ such that $e''$ is contained in $a$ sets of $\mathcal{F}'$ that must still be contained in $a$ sets of $\mathcal{F}''$. If not, each element of the new family $\mathcal{F}''$ would be contained in at most $a-1$ sets of $\mathcal{F}''$, and so we would have that $f(n,a-1)\geq m+k-k'=m$, which is a contradiction on the fact that $f(n,a-1) < f(n,a)$.

Thus if $k=k'$, then there exists $e''$ such that $|\{S\in \mathcal{F}'| S\ni e''\}|=a$ and such that $e''$ is never contained in a set of the family that does not contain $e'$. Indeed, if there existed $S'\in \mathcal{F}'$ such that $e'\not\in S'$ and $e''\in S'$, then for any set $S''\in \mathcal{F}$ such that $S''\cup e' \in \mathcal{F}$ as well, i.e. twin sets with difference $e'$, then one of $S'\cup S''$ and $S'\cup (S''\cup e')$ will disappear in $\mathcal{F}''$ and so $e''$ would be present $a-1$ times, a contradiction.

Thus $\{S\in \mathcal{F}'| S\ni e''\} \subseteq \{S\in \mathcal{F}'| S\ni e'\}$ and since $|\{S\in \mathcal{F}'| S\ni e''\}|=a$, then $|\{S\in \mathcal{F}'| S\ni e'\}|=a$ as well and $\{S\in \mathcal{F}'| S\ni e''\} = \{S\in \mathcal{F}'| S\ni e'\}$ which is a contradiction of the fact that $e'$ and $e''$ are twin differences for some sets. Since they are copies of each other, we can remove either one of them without creating duplicate sets.
Thus, $k < k'$, and so if $k'=1$, then $k\leq 0$, and then $f(n,a)\geq f(n+1,a)$.
\qed

Note that this also means that if we remove any element from such a solution, and remove a copy of every duplicated set created, what remains is never an optimal solution for $f(n,a)$.

\section{Conclusion}

As we have seen, a complete proof of the new conjectures has far-reaching implications: the Frankl conjecture would hold for about $\frac{2}{3}$ of all possible cases, and we could show that there always exists an element in $\frac{6}{13}$ of the sets of a union-closed family. Therefore we believe our new conjectures merit additional attention. In order to encourage further progress in this direction, we conclude with a few open problems of interest.

\begin{enumerate}
\item Show that $f(n,a)=f(n+1,a)$ for smaller values of $n$, i.e., for values of $n<a$.
\item Show that if $f(n,a)<f(n+1,a)$, each element is a twin difference for even more sets. As noted, at some point, this clearly implies that the solutions is not optimal or even feasible. 
\item Find a constant upper bound for $f_t(n+1,a)-f_t(n,a)$. Since we know $f(n,a)$ stops growing after $n=a$, this would provide a first constant lower bound for the number of sets containing the most frequent element in a union-closed family.
\end{enumerate}
\newpage
\input{appendix} \newpage
\nocite{*}
\bibliography{franklrevised}
\bibliographystyle{alpha}

\end{document}

%% file: appendix.tex
\section{Appendix}

In the following tables for $f(n,a)$ and $g(n,m)$, we remove unnecessary columns, i.e. the columns for which the values of $f(n,a)$ and $g(n,m)$ are trivial.

\begin{table}[h!]
 \caption{Values of $f(n,a)$}
 \label{tablef}
\begin{center}
\begin{scriptsize}
\setlength{\tabcolsep}{2pt}
\begin{tabular}{|l|ccccccccc|}
\hline
$a\backslash n$ & 1 & 2 & 3 & 4 & 5& 6& 7& 8 & 9\\
  \hline
  1 & 2 & 2 & 2 & 2 & 2& 2& 2& 2 & 2 \\
  2 &  & 4 & 4 & 4& 4& 4& 4& 4 & 4\\
  3 & & & 5 & 5 & 5 & 5 & 5 & 5 & 5\\
  4 & & & 8& 8& 8& 8& 8& 8 & 8\\ 
  5 & & & & 9& 9& 9& 9& 9 & 9\\
  6 & & & & 10& 10& 10& 10& 10 & 10\\
  7 & & & & 12& 12& 12& 12& 12 & 12\\
  8 & & & & 16& 16& 16& 16& 16 & 16\\
  9 & & & & & 17& 17& 17& 17 & 17\\
  10 & & & & & 18& 18& 18& 18 & 18\\
  11 & & & & & 19& 19& 19& 19 & 19 \\
  12 & & & & & 21& 21& 21& 21 & 21\\
  13 & & & & & 23& 23& 23& 23 & 23\\
  14 & & & & & 25 & 25 & 25 & 25 & 25\\
  15 & & & & & 27& 27& 27& 27 & 27\\
  16 & & & & & 32& 32& 32& 32 & 32\\
  17 & & & & & & 33 & 33 & 33 & 33\\
  18 & & & & & & 34 & 34 & 34 & 34\\
  19 & & & & & & 35 & 35 & 35 & 35\\
  20 & & & & & & 36 & 36 & 36 & 36\\
  21 & & & & & & 38 & 38 & 38 & 38\\
  22 & & & & & & 40 & 40 & 40 & 40\\
  23 & & & & & & 41 & 41 & 41 & 41\\
  24 & & & & & & 43 & 42 & 42 & 42\\
  25 & & & & & & 45 & 45 & 45 & 45\\
  26 & & & & & & 47 & 47 & 47 & 47\\
  27 & & & & & & 49 & 49 & 49 & 49\\
  28 & & & & & & 52 & 52 & 52 & 52\\
  29 & & & & & & 53 & 53 & 53 & 53\\
  30 & & & & & & 56 & 56 & 56 & 56\\
  31 & & & & & & 58 & 58 & 58 & 58\\
  32 & & & & & & 64 & 64 & 64 & 64\\
  \hline
\end{tabular} \quad \begin{tabular}{|l|ccc|}
\hline
$a\backslash n$ & 7& 8 & 9\\
  \hline
33 & 65& 65& 65\\
34 & 66& 66& 66\\
35 & 67& 67& 67\\
36 & 68& 68& 68\\
37 & 69& 69& 69\\
38 & 71& 71& 71\\
39 & 72& 72& 72\\
40 & 74& 74& 74\\
41 & 75& 75& 75\\
42 & 77& 77& 77\\
43 & 79& 79& 79\\
44 & 80& 80& 80\\
45 & 82& 82& 82\\
46 & 83& 83& 83\\
47 & 85& 85& 85\\
48 & 88& 88& 88\\
49 & 89& 89& 89\\
50 & 91& 91& 91\\
51 & 93& 93& 93\\
52 & 95& 95& 95\\
53 & 98& 98& 98 \\
54 & 99& 99& 99\\
55 & 101& 101& 101\\
56 & 104& 104& 104\\
57 & 105& 105& 105\\
58 & 108& 108& 108\\
59 & 110& 110& 110\\
60 & 113& 113& 113\\
61 & 115& 115& 115\\
62 & 118& 118& 118\\
63 & 121& 121& 121\\
64 & 128& 128& 128\\
  \hline
\end{tabular} \quad \begin{tabular}{|l|cc|}
\hline
$a\backslash n$ & 8 & 9\\
  \hline
65 & 129& 129\\
66 & 130& 130\\
67 & 131& 131\\
68 & 132& 132\\
69 & 133& 133\\
70 & 134& 134\\
71 & 136& 136\\
72 & 137& 137\\
73 & 139& 139\\
74 & 140& 140\\
75 & 142& 142\\
76 & 144& 144\\
77 & 145& 145\\
78 & 146& 146\\
79 & 147& 147\\
80 & 149& 149\\
81 & 150& 150\\
82 & 152& 152\\
83 & 154& 154\\
84 & 156& 156\\
85 & 157& 157\\
86 & 158& 158\\
87 & 160& 160\\
88 & 162& 162\\
89 & 164& 164\\
90 & 166& 166\\
91 & 168& 168\\
92 & 170& 170\\
93 & 171& 171\\
94 & 173& 173\\
95 & 175& 175\\
96 & 176& 176\\
  \hline
\end{tabular} \quad \begin{tabular}{|l|cc|}
\hline
$a\backslash n$ & 8 & 9\\
  \hline
97 & 179& 179\\
98 & 180& 180\\
99 & 182& 182\\
100 & 184& 184\\
101 & 186& 186\\
102 & 188& 188\\
103 & 189& 189\\
104 & 192& 192\\
105 & 194& 194\\
106 & 196& 196\\
107 & 198& 198\\
108 & 200& 200\\
109 & 202& 202\\
110 & 204& 204\\
111 & 206& 206\\
112 & 209& 209\\
113 & 211& 211\\
114 & 214& 214\\
115 & 216& 216\\
116 & 220& 220\\
117 & 221& 221\\
118 & 224& 224\\
119 & 226& 226\\
120 & 229& 229\\
121 & 231& 231\\
122 & 233& 233\\
123 & 236& 236\\
124 & 240& 240\\
125 & 242& 242\\
126 & 245& 245\\
127 & 248& 248\\
128 & 256& 256\\
  \hline
\end{tabular}
\end{scriptsize}
\end{center}
\end{table}

\newpage 
\bigskip
\bigskip
\begin{table}[h!]
\caption{Values of $g(n,m)$}
\label{tableg}
\begin{center}
\begin{footnotesize}
\setlength{\tabcolsep}{2pt}
\begin{tabular}{|l|cccccccc|}
\hline
$m\backslash n$  & 1 & 2 & 3 & 4 & 5& 6& 7& 8 \\
  \hline
   &  &  &  &  & & & &   \\
  2 &  1& 1 &1 & 1& 1& 1& 1& 1 \\
  3 & & 2& 2 & 2 & 2 & 2 & 2 & 2 \\
  4 & & 2& 2& 2& 2& 2& 2& 2\\ 
  5 & & &3 & 3& 3& 3& 3& 3\\
  6 & & & 4& 4& 4& 4& 4& 4 \\
  7 & & & 4& 4&4& 4& 4& 4\\
  8 & & & 4& 4& 4& 4& 4& 4\\
  9 & & & & 5& 5& 5& 5& 5\\
  10 & & & & 6& 6& 6& 6& 6\\
  11 & & & & 7& 7& 7& 7& 7\\
  12 & & & & 7& 7& 7& 7& 7\\
  13 & & & & 8& 8& 8& 8& 8\\
  14 & & & & 8& 8 & 8 & 8 & 8\\
  15 & & & & 8& 8& 8& 8& 8\\
  16 & & & & 8& 8& 8& 8& 8\\
17 & & & & & 9& 9& 9& 9\\
18 & & & & & 10& 10&10&10\\
19 & & & & & 11& 11& 11& 11\\
20 & & & & & 12& 12& 12& 12\\
21 & & & & & 12& 12& 12& 12\\
22 & & & & & 13& 13& 13& 13\\
23 & & & & & 13& 13& 13& 13\\
24 & & & & & 14& 14& 14& 14\\
25 & & & & & 14& 14& 14& 14\\
26 & & & & & 15& 15& 15& 15\\
27 & & & & & 15& 15& 15& 15\\
28 & & & & & 16& 16& 16& 16\\
29 & & & & & 16& 16& 16& 16\\
30 & & & & & 16& 16& 16& 16\\
31 & & & & & 16& 16& 16& 16\\
32 & & & & & 16& 16& 16& 16\\
  \hline
\end{tabular} \quad \begin{tabular}{|l|ccc|}
\hline
$m\backslash n$  & 6& 7& 8 \\
  \hline
33 & 17&17 &17 \\
34 & 18& 18& 18\\
35 & 19& 19& 19\\
36 & 20& 20& 20\\
37 & 21& 21& 21\\
38 & 21& 21& 21\\
39 & 22& 22& 22\\
40 & 22& 22& 22\\
41 & 23&23 & 23\\
42 & 24&24 &24 \\
43 & 24& 24& 24\\
44 & 25& 25& 25\\
45 & 25& 25& 25\\
46 & 26& 26& 26\\
47 & 26& 26& 26\\
48 & 27& 27& 27\\
49 & 27& 27& 27\\
50 & 28& 28& 28\\
51 & 28& 28& 28\\
52 & 28& 28& 28\\
53 & 29& 29& 29\\
54 & 30& 30& 30\\
55 & 30& 30& 30\\
56 & 30& 30& 30\\
57 & 31& 31& 31\\
58 & 31& 31& 31\\
59 & 32& 32& 32\\
60 & 32& 32& 32\\
61 & 32& 32& 32\\
62 & 32& 32& 32\\
63 & 32& 32& 32\\
64 & 32& 32& 32\\
  \hline
\end{tabular} \quad \begin{tabular}{|l|cc|}
\hline
$m\backslash n$ & 7& 8 \\
  \hline
65 & 33& 33\\
66 & 34& 34\\
67 & 35& 35\\
68 & 36& 36\\
69 & 37& 37\\
70 & 38& 38\\
71 & 38& 38\\
72 & 39& 39\\
73 & 40& 40\\
74 & 40& 40\\
75 & 41& 41\\
76 & 42& 42\\
77 & 42& 42\\
78 & 43& 43\\
79 & 43& 43\\
80 & 44& 44\\
81 & 45& 45\\
82 & 45& 45\\
83 & 46& 46\\
84 & 47& 47\\
85 & 47& 47\\
86 & 48& 48\\
87 & 48& 48\\
88 & 48& 48\\
89 & 49& 49\\
90 & 50& 50\\
91 & 50& 50\\
92 & 51& 51\\
93 & 51& 51\\
94 & 52& 52\\
95 & 52& 52\\
96 & 53& 53\\
  \hline
\end{tabular} \quad \begin{tabular}{|l|cc|}
\hline
$m\backslash n$ & 7& 8 \\
  \hline
97 & 53& 53\\
98 & 53& 53\\
99 & 54& 54\\
100 & 55& 55\\
101 & 55& 55\\
102 & 56& 56\\
103 & 56& 56\\
104 & 56& 56\\
105 & 57& 57\\
106 & 58& 58\\
107 & 58& 58\\
108 & 58& 58\\
109 & 59& 59\\
110 & 59& 59\\
111 & 60& 60\\
112 & 60& 60\\
113 & 60& 60\\
114 & 61& 61\\
115 & 61& 61\\
116 & 62& 62\\
117 & 62& 62\\
118 & 62& 62\\
119 & 63& 63\\
120 & 63& 63\\
121 & 63& 63\\
122 & 64& 64\\
123 & 64& 64\\
124 & 64& 64\\
125 & 64& 64\\
126 & 64& 64\\
127 & 64& 64\\
128 & 64& 64\\
  \hline
\end{tabular}
\end{footnotesize}
\end{center}
\end{table}